\documentclass{article}

\usepackage{fullpage} % Page layout

\usepackage{fancyhdr} % Headings
\usepackage{manfnt} % Dangerous bend symbols \dbend \textdbend \lhdbend \textlhdbend 
\usepackage{latexsym,stmaryrd,amsfonts,amssymb,amsmath,amsthm,amsbsy,array} % Mathematics symbols
\usepackage{graphicx} % Figures
\usepackage{amscd,xy} % Diagrams
\usepackage{tikz-cd} 
\usepackage{mathtools}

\usepackage[OT2,OT1]{fontenc} % Cyrillic script
\newcommand\cyr{%
\renewcommand\rmdefault{wncyr}%
\renewcommand\sfdefault{wncyss}%
\renewcommand\encodingdefault{OT2}%
\normalfont
\selectfont}
\DeclareTextFontCommand{\textcyr}{\cyr}

\def \what{\widehat}

\def \Z{\mathbb{Z}}
\def \Q{\mathbb{Q}}

\def \A{\mathcal{A}}

\def \Cg1{\mathcal{C}_{g,1}}

\def \CSg1{\mathcal{C}(\Sigma_{g,1})}

\def \ICSg1{\mathcal{IC}(\Sigma_{g,1})}

\begin{document}

\title{Jacobi diagrams on surfaces and quantum invariants}
\date{}
\author{Dorin Cheptea\thanks{This work was partially supported by 
a grant of the Romanian National Authority for Scientific Research, 
CNCS-UEFISCDI, project number PN-II-RU-TE-2012-3-0492.}}

\maketitle

\begin{abstract}

This is a preliminary version; it will be completed shortly.


\end{abstract}

\section{Intorduction}

The relationship between quantum link invariants (which generalize the Jones polynomial),
and the Kontsevich integral (universal finite type invariant of knots/links) is fairly well understood.
The quantum link invariants were extended to 3-manifolds by Witten and Reshetikhin-Turaev (WRT). 
One associates to a compact oriented 3-manifold $M$, a root of unity $\xi$, 
and a semi-simple Lie algebra $\mathfrak g$, for example ${\mathfrak sl}_2$,
a complex number $\tau_M(\xi)$.

The Kontsevich integral was extended to an invariant of 3-manifolds (LMO) by Le-Murakami-Ohtsuki.
LMO takes values in a certain algebra $\A(\emptyset)$ of Jacobi diagrams.
Any finite-dimensional semi-simple metrized Lie algebra $\mathfrak g$
defines a linear map $\A(\emptyset)\rightarrow\Q[[h]]$ 
called the weight system associated to $\mathfrak g$. 
For a Jacobi diagram $D\in\A(\emptyset)$, 
$\what{W}_{{\mathfrak sl}_2}(D)=\sum_{d\geq 0}W_{{\mathfrak sl}_2}(D)h^d$,
where $W_{{\mathfrak sl}_2}(D)\in\Q$, $d=$ degree of $D$.

The relationship between WRT and LMO invariants is more complex, 
and is best described by the theory \cite{Habiro2008} of the unified invariant $J_M(q)$,
whose evaluation at any root of unity $q=\xi$ conincides 
with the value of the WRT invariant at that root: $J_M(\xi)=\tau_M(\xi)$.
Habiro's invariant takes values in the Habiro's ring:
$$\what{\Z[q]}:=\lim\limits_{\leftarrow n}\frac{\Z[q]}{((1-q)(1-q^2)\cdots(1-q^n))},$$
the cyclotomic completion of $\Z[q]$.
Every element $f(q)\in\what{\Z[q]}$ can be written (not uniquely) as an infinite sum
$f(q)=\sum_{k\geq 0}f_k(q)(1-q)(1-q^2)\cdots(1-q^k)$,
with $f_k(q)\in\Z[q]$.
When $q=\xi$ a root of unity, only a finite number of these terms are not zero, hence the evaluation
$ev_{\xi}(f(q))$ is well-defined. Moreover, it only depends on $f(q)$.
Habiro's ring has the property that the formal Taylor series expansion of $f(q)\in\what{\Z[q]}$ 
at a root $\xi$ of $1$
$$T_{\xi}(f)=\sum\limits_{n=0}^{\infty}\frac{f^(n)(\xi)}{n!}(q-\xi)^n$$
is well defined. Moreover, the map $T_{\xi}:\what{\Z[q]}\rightarrow\Z[\xi][[q-\xi]]$ is injecctive,
i.e. a function in $\what{\Z[q]}$ is determined by its Taylor expansion at a point in the domain $U$,
the set of roots of $1$.
While the Taylor series $T_1f$ has convergence radius zero, 
in the $p$-adic topology $T_1f(\xi)$ converges to $f(\xi)$. 

Properties of the Habiro's ring imply that the following diagram is commutative \cite{Habiro2008,BL2011}:
$$\begin{tikzcd}
\{\Z HS\} \arrow{rd}[swap]{LMO} \arrow{r}{J_M(q)} & \what{\Z[q]} \arrow[hook]{r}{T_1} 
& \Z[[1-q]] \arrow[hook]{d}{h=1-q} \\
& \A(\emptyset) \arrow{r}{{\mathfrak sl}_2-weight}[swap]{system} & \Q[[h]] 
\end{tikzcd}$$
where $\{\Z HS\}$ denotes the set of integral homology 3-spheres. This puts into perspective the
well-known Ohtsuki's result that for an integer homology 3-sphere $M$, the perturbative
$PSU(2)=SO(3)$ invariant (see \cite{Ohtsuki2002}) recovers from the LMO invariant as
$\tau^{SO(3)}(M)=\what{W}_{{\mathfrak sl}_2}(Z^{LMO}(M))$,
while also proving that $\tau^{SO(3)}(M)$ has integer coefficients
(a priori $\tau^{SO(3)}(M)\in\Q[[q-1]]$).

Both WRT and LMO invariants extend to 3-manifolds with boundary in a nice way (TQFTs).
However, the two types of functorialities have different flavor.
The LMO functor is appropriate for homology cylinders, which are a generalization to cobordisms
of integer homology 3-spheres. The TQFTs fot the WRT invariants, by contrast,
make use of the totality of 3-manifolds with boundary.

The weight system construction generalizes to a map
$W_{\mathfrak g}:(\A(H_{\Q}),\star)\rightarrow(S({\mathfrak g}\otimes H_{\Q})[[t]],\star)$
from the Hopf algebra of symplectic Jacobi diagrams,
where $H_{\Q}=H\otimes\Q$, and $H=H_1(\Sigma_{g,1},\Z)$,
to the deformation qunatization of the symmetric algebra $S({\mathfrak g}\otimes H_{\Q})$,
whose Poisson structure is induced by the symplectic form on $H_{\Q}$ (see \cite{HM2009}).
The LMO homomorphism sends the monoid of homology cylinders $\ICSg1$ (see Section 2 below)
to the group-like elements $GLike\:\A(H_{\Q})$,
and composition of cobordisms is sent to the multiplication $\star$ which can be described in
purely algebraic-combinatorial terms (see \cite{HM2012}). (On a side note, the induced
Lie bracket $[,]_{\star}$ on the reduction of $\A(H_{\Q})$ to connected and tree-like
Jacobi diagrams $\A^{t,c}(H_{\Q})$ is the one considered independently by Kontsevich and Morita in 1993.)
The map $W_{\mathfrak g}$ sends $\star$ to the Moyal-Weyl product on 
$S({\mathfrak g}\otimes H_{\Q})[[t]]$ (see \cite{HM2009}).

This provides a means by which some functoriality properties of WRT and LMO invariants can be compared.
The aim of this paper is to introduce (adapt) another ingredient to this goal:
Jacobi diagrams on surfaces, which have been considered by Andersen-Mattes-Reshetikhin \cite{AMR1998}.

The paper is organized as follows. Section 2 is an overview of typies of cobordisms over the surfaces 
$\Sigma_{g,b}$ that we will be considering.
Jacobi diagrams on surfaces are addressed in Section 3.
In Section 4 we introduce a new map modeled on the LMO homomorphism.
Section 5 ...

\section{Cobordisms over $\Sigma_{g,b}$}

Let $\Sigma_{g,b}$ denote a compact connected oriented surface of genus $g$ with $b$ boundary components.
For our purposes, $b$ will be $0$ or $1$.
Let $H_g$ denote a handlebody of genus $g$ such that $\partial H_g=\Sigma_g$.
The mapping class group of $\Sigma_{g,b}$ is the group of isotopy classes of homeomorphisms
$\Sigma_{g,b}\rightarrow\Sigma_{g,b}$ which preserve the orientation and fix the boundary pointwise.
It is well-known that every compact connected oriented 3-manifold is homeomorphic to $H_g\cup_fH_g$
for some orientation preserving homeomorphism $f$ of $\Sigma_b$.
The Torelli group of $\Sigma_{b,g}$ is the subgroup of the mapping class group, whose elements are
represented by homeomorphisms that induce the identity map on the first homology $H_1(\Sigma_{g,b})$.
It is well-known that every (compact connected oriented) integer homology 3-sphere is homeomorphic 
to $H_g\cup_fH_g$ for some orientation preserving homeomorphism $f$ of $\Sigma_b$, who isotopy class is
in the Torelli group.

We shall call a pair $(M,m)$ a cobordism of $\Sigma_{g,b}$ if
$M$ is a compact connected oriented 3-manifold and 
$m:\partial(\Sigma_{g,b}\times[-1,1])\rightarrow\partial M$
is an orientation-preserving homeomorphism.
$(M,m)$ and $(M^{\prime},m^{\prime})$ are equivalent if there is an orientation-preserving homeomorphism
$f:M\rightarrow M^{\prime}$ such that $f|_{\partial M}\circ m=m^{\prime}$.
$m_{\pm}=m(\cdot,\pm1):\Sigma_{g,b}\rightarrow M$ determines {\it top} and {\it bottom} surfaces, 
and are used to compose cobordisms:
$M\circ M^{\prime}=M\cup_{m_+\circ(m^{\prime}_-)^{-1}}M^{\prime}$.
Here $(\Sigma_{g,b}\times[-1,1], id)$ is the identity element of resulting monoid $Cob(Sigma_{g,b})$.

Notice the submonoids of homology cobordisms and homology cylingers over $\Sigma_{g,b}$, 
where a cobordism $(M,m)$ is called a homology cobordism if both $m_+$ and $m_-$ induce isomprphisms 
$H_{\ast}(Sigma_{g,b})\rightarrow H_{\ast}(M)$, and a homology cylinder if in addition the two induced
isomorphisms coincide.

The mapping cylinder construction associates to a surface homeomorphisms $f$ the cobordism
$(\Sigma_{g,b}\times[-1,1],id\times\{-1\}\cup f\times\{1\})$.
Thus, Torelli group is contained in the monoid of homology cylinders.

\section{Jacobi diagrams on surfaces}

\section{The LMO map}

\section{}

\bigskip
\noindent
{\sc Institute of Mathematics of the Romanian Academy,\\ P. O. Box 1-764, Bucharest, 014700, Romania}

\noindent
{\it Email}: Dorin.Cheptea@imar.ro

\end{document}